# AN ADAPTIVE STEP-DOWN PROCEDURE WITH PROVEN FDR CONTROL UNDER INDEPENDENCE

BY YULIA GAVRILOV,[1] YOAV BENJAMINI[1] AND SANAT K. SARKAR[2]

*Tel Aviv University, Tel Aviv University and Temple University*

In this work we study an adaptive step-down procedure for testing $m$ hypotheses. It stems from the repeated use of the false discovery rate controlling the linear step-up procedure (sometimes called BH), and makes use of the critical constants $iq/[(m+1-i(1-q)]$, $i=1,\dots,m$. Motivated by its success as a model selection procedure, as well as by its asymptotic optimality, we are interested in its false discovery rate (FDR) controlling properties for a finite number of hypotheses. We prove this step-down procedure controls the FDR at level $q$ for independent test statistics. We then numerically compare it with two other procedures with proven FDR control under independence, both in terms of power under independence and FDR control under positive dependence.

**1. Introduction.** Multiple testing has become again a flourishing area of research. With new challenges stemming from very large and complex applied problems, with implications to seemingly unrelated areas such as model selection and signal denoising, it is no wonder that a large variety of approaches and methods, sometimes competing and other times complementing, are being developed.

Within the family of methods whose concern is about the false discovery rate (FDR), we follow here a rudimentary suggestion in Benjamini, Krieger and Yekutieli [4] (hereafter referred to as BKY), for testing $m$ null hypotheses $H_i$, $i=1,\dots,m$, using a $p$-value $p_i$ available for each $H_i$. We make use of the constants

$$\alpha_i = iq/(m+1-i(1-q)), \qquad i=1,\dots,m, \tag{1.1}$$

in a step-down manner.

Received June 2007; revised December 2007.
[1]Supported by G.I.F. Grant 847/2004 and by NIH Grant DA015087.
[2]Supported by NSF Grant DMS-06-03868.

More formally, we define the multiple-stage adaptive step-down procedure as follows: after sorting the $p$-values as $p_{(1)} \leq \cdots \leq p_{(m)}$, with $H_{(i)}$ being the null hypothesis corresponding to $p_{(i)}$, we let

$$k = \max\{1 \leq i \leq m : p_{(j)} \leq \alpha_j, j = 1, \ldots, i\}$$

and reject the $k$ hypotheses associated with $p_{(1)}, \ldots, p_{(k)}$, when such a $k$ exists; otherwise, reject no hypothesis.

Motivated by the success of this procedure as a model selection procedure [2], as well as by its asymptotic optimality in a mixture model noted in [7], we consider investigating its FDR controlling properties. To that end, we first theoretically establish its FDR control under independence by proving the following slightly more general theorem.

THEOREM 1.1. *For independent test statistics, with the $p$-values distributed as $U(0,1)$ under the corresponding null hypotheses, the step-down procedure with the critical values $\alpha_1 \leq \cdots \leq \alpha_m$ satisfying $\alpha_i/(1-\alpha_i) \leq iq/(m+1-i)$, $i = 1, \ldots, m$, controls the FDR at $q$.*

We then numerically compare its power to other compatible procedures with proven FDR control under independence. We similarly examine its FDR controlling behavior under positive dependence, and find it to control the FDR mostly below the desired level, but sometimes slightly above it. The possibility of choosing different sets of constants that satisfy the condition of Theorem 1.1 and further offer proven FDR control under positive dependence is discussed. One such procedure is suggested, but its proven FDR under positive dependency comes at the expense of being too conservative.

## 2. Background.

This work is within the FDR framework, where the interest is in procedures that control the expected proportion of true null hypotheses rejected ($V$) out of the total number of null hypotheses rejected ($R$), in the form of $\text{FDR} = E[(V/R)I_{\{R>0\}}]$. To understand the rationale behind the development of an FDR controlling adaptive procedure, recall first that when testing $m$ hypotheses, the linear step-up procedure at level $q$ in Benjamini and Hochberg [3], also called the BH procedure, runs as follows: let $H_{(i)}$ correspond to the ordered $p$-value $p_{(i)}$, as before, and

$$k = \max\{1 \leq i \leq m : p_{(i)} \leq iq/m\}.$$

If such a $k$ exists, reject the $k$ hypotheses associated with $p_{(1)}, \ldots, p_{(k)}$; otherwise reject none. When $m_0$ of the hypotheses are true, the BH procedure controls the FDR at $qm_0/m$ when the test statistics are independent. It is too conservative when $m_0/m$ is small. Knowledge of $m_0$ can therefore be very useful to improve upon the performance of the BH procedure. If $m_0$ was given to us by an “oracle,” the linear step-up procedure with $q' = qm/m_0$ would control the FDR at precisely the desired level $q$. In practice $m_0$ is not known, and so the field is open for efforts to construct *adaptive procedures* that estimate this factor first and make use of it in the BH procedure.

In BKY, a two-stage adaptive BH procedure is given based on this rationale. It estimates $m_0$ at the first stage based on the number of null hypotheses rejected by the BH procedure, then uses this estimate to modify the BH procedure at the second stage. They prove that this two-stage BH procedure controls the FDR under independence. They also offer two multiple-stage versions of this, where the BH procedure is used in a similar way up to $m$ times rather than merely twice. The first procedure runs as follows: let

$$k = \max\{1 \le i \le m : \forall j \le i, \exists\, l \ge j \text{ so that } p_{(l)} \ge lq/(m+1-j(1-q))\}.$$

If such a $k$ exists, reject the $k$ hypotheses associated with $p_{(1)}, \ldots, p_{(k)}$; otherwise reject no hypothesis. A simpler version of the above procedure requires $l = j$, resulting in a multiple-stage step-down procedure with the constants $\{\alpha_i\}$ in (1.1). Both of these procedures compare the ordered $p$-values to the same constants, and BKY offer no analytical proof of the FDR controlling properties of either procedure. In the first one, a combination of step-up and step-down methods is used. In the second simpler one, a pure step-down method is used. It is this second multiple-stage adaptive step-down procedure that we revisit in this article, offering an analytical proof of its FDR control under independence. It is interesting to note that the FDR is not controlled by the corresponding step-up procedure based on these constants. As the last constant in (1.1) approaches 1 as $m$ increases, almost all tested hypotheses will be rejected at the first stage in this step-up procedure, and thus the FDR will be approximately $m_0/m$.

Following the work of Abramovich, Benjamini, Donoho and Johnstone [1], who showed the asymptotic advantage of using BH-based penalty functions for model selection in sparse situations, Benjamini and Gavrilov [2] investigated via finite sample simulations the performance of model selection procedures based on the FDR control over a wide range of configurations. Their analysis revealed excellent performance for the model selection method based on the constants $\{\alpha_i\}$ in (1.1). Moreover, the method can be implemented as a forward selection penalized least squares, just as the BH procedure can be implemented within a backward elimination procedure. In view of these results Benjamini and Gavrilov [2] were able to prove the FDR controlling properties of the multiple-stage step-down procedure, but only in an asymptotic setting.

Independently, within the same asymptotic framework of a mixture model, Finner, Dickhaus and Roters [7] were interested in the form of an asymptotically optimal curve to which the empirical distribution of the $p$-values should be compared. They arrived at a curve of constants to which the inverse of the empirical distribution is to be compared that asymptotically takes the form $f(t) = tq/[1 - t(1 - q)]$. As $\alpha_i = f(i/(m+1))$, the equivalence of the curve in [7] and the curve to which the constants of the multiple-stage adaptive step-down procedure converge as $m$ increases, is clear. Moreover, when searching for a finite procedure, they looked at methods that use the constants $\alpha_{i,m} = iq/(m + \beta_m - i(1 - q))$, with $\beta_m/m \to 0$. The $\alpha_i$ in (1.1) corresponds to $\beta_m = 1$ in these constants. They were also able to prove that a few variations of this procedure control asymptotically the FDR. They nevertheless noted that getting a proof for the finite case of the step-down procedure that uses these constants seems to be a very difficult task.

This brings us to the result presented in this paper in which we prove that for any finite number of hypotheses, tested using independent test statistics, the multiple-stage adaptive step-down procedure controls the FDR.

The present adaptive procedure should be compared with two other adaptive procedures, as they too have been proven to offer FDR control when the test statistics are independent. These are the two-stage adaptive BH procedure in [4] and the modified version of Storey's procedure [9] presented in [10]. The two-stage BH procedure runs formally as follows: at the first stage apply the linear step-up procedure at level $q' = q/(1 + q)$. Let $r_1$ be the number of rejected hypotheses. If $r_1 = 0$ reject no hypotheses and stop; if $r_1 = m$ reject all $m$ hypotheses and stop. Otherwise, let $\hat{m}_0 = m - r_1$, and use the linear step-up procedure at the second stage with $q^* = q'm/\hat{m}_0$. The modified Storey's procedure runs as follows: let $r(\lambda) = \sharp\{p_i \le \lambda\}$, then estimate $m_0$ by $\hat{m}_0 = (m + 1 - r(\lambda))/(1 - \lambda)$, and use the linear step-up procedure with $q^* = qm/\hat{m}_0$ with further requiring for selection that $p_i \le \lambda$. For estimating $m_0$, Storey [9] and Storey and Tibshirani [11] recommended use of $\lambda = 0.5$.

**3. Proof of Theorem 1.1.** In this section, we provide a proof of Theorem 1.1. Let us now introduce the more explicit notation $P_{i:m} = P_{(i)}$, so $P_{1:n} \le \cdots \le P_{m:m}$ are the $m$ ordered random $p$-values, and $\{H_i, i \in I_0\}$ is the set of true null hypotheses. Then, from Sarkar [8] we see that the FDR of the step-down procedure in the theorem is given by

$$\begin{aligned}
\mathrm{FDR} &= \sum_{i \in I_0} \sum_{r=1}^{m} \frac{1}{r} \Pr\{P_{1:m} \le \alpha_1, \ldots, P_{r:m} \le \alpha_r, P_{r+1:m} \ge \alpha_{r+1}, P_i \le \alpha_r\} \\
&= \sum_{i \in I_0} \sum_{r=1}^{m} E\biggl[\Pr\{P_{1:m} \le \alpha_1, \ldots, P_{r:m} \le \alpha_r \mid P_i\} \\
&\qquad \cdot \left\{\frac{I(P_i \le \alpha_r)}{r} - \frac{I(P_i \le \alpha_{r-1})}{r-1}\right\}\biggr] \qquad (3.1)
\end{aligned}$$

$$= \sum_{i \in I_0} \sum_{r=1}^{m} \frac{1}{r} \Pr\{P_{1:m} \le \alpha_1, \ldots, P_{r:m} \le \alpha_r, \alpha_{r-1} < P_i \le \alpha_r\}$$
$$- \sum_{i \in I_0} \sum_{r=2}^{m} \frac{1}{r(r-1)} \Pr\{P_{1:m} \le \alpha_1, \ldots, P_{r:m} \le \alpha_r, P_i \le \alpha_{r-1}\}$$

(with $P_{m+1:m} = 1$, $\alpha_0 = 0$, $\alpha_{m+1} = 1$ and $\alpha_0/0 = 0$).

Let $P^{(-i)}_{1:m-1} \le \cdots \le P^{(-i)}_{m-1:m-1}$ be the ordered versions of the $m-1$ $p$-values in the set $(P_1, \ldots, P_m) \setminus \{P_i\}$. Using the following result:

$$\begin{aligned} &\{P_{1:m} \le c_1, \ldots, P_{m:m} \le c_m, c_{j-1} < P_i \le c_j\} \\ &\quad = \{P^{(-i)}_{1:m-1} \le c_1, \ldots, P^{(-i)}_{j-1:m-1} \le c_{j-1}, P^{(-i)}_{j:m-1} \le c_{j+1}, \ldots, \\ &\qquad P^{(-i)}_{m-1:m-1} \le c_m, c_{j-1} < P_i \le c_j\}, \end{aligned} \tag{3.2}$$

$j = 1, \ldots, m$, for any increasing set of constants $c_0 = 0 \le c_1 \le \cdots \le c_m \le 1$ (see, e.g., [6]), we notice that

$$\begin{aligned} &\Pr\{P_{1:m} \le \alpha_1, \ldots, P_{r:m} \le \alpha_r, \alpha_{r-1} < P_i \le \alpha_r\} \\ &\quad = \Pr\{P^{(-i)}_{1:m-1} \le \alpha_1, \ldots, P^{(-i)}_{r-1:m-1} \le \alpha_{r-1}, \alpha_{r-1} < P_i \le \alpha_r\} \end{aligned} \tag{3.3}$$

for $1 \le r \le m$ [with $P^{(-i)}_{0:m-1} = 0$] and

$$\begin{aligned} &\Pr\{P_{1:m} \le \alpha_1, \ldots, P_{r:m} \le \alpha_r, P_i \le \alpha_{r-1}\} \\ &\quad \ge \Pr\{P^{(-i)}_{1:m-1} \le \alpha_1, \ldots, P^{(-i)}_{r-1:m-1} \le \alpha_{r-1}, P_i \le \alpha_{r-1}\} \end{aligned} \tag{3.4}$$

for $2 \le r \le m$.

Thus, using (3.3) and (3.4) in (3.1), we have

$$\begin{aligned} \mathrm{FDR} &\le \sum_{i \in I_0} \sum_{r=1}^{m} \Pr\{P^{(-i)}_{1:m-1} \le \alpha_1, \ldots, P^{(-i)}_{r-1:m-1} \le \alpha_{r-1}\} \left\{ \frac{\alpha_r}{r} - \frac{\alpha_{r-1}}{r-1} \right\} \\ &= \sum_{i \in I_0} \sum_{r=1}^{m} \frac{\alpha_r}{r} \Pr\{P^{(-i)}_{1:m-1} \le \alpha_1, \ldots, P^{(-i)}_{r-1:m-1} \le \alpha_{r-1}, P^{(-i)}_{r:m-1} > \alpha_r\} \\ &\le \alpha \sum_{i \in I_0} \sum_{r=1}^{m} \frac{1-\alpha_r}{m-r+1} \Pr\{P^{(-i)}_{1:m-1} \le \alpha_1, \ldots, P^{(-i)}_{r-1:m-1} \le \alpha_{r-1}, \\ &\qquad P^{(-i)}_{r:m-1} > \alpha_r\} \\ &= \alpha \sum_{i \in I_0} \sum_{r=1}^{m} \frac{1}{m-r+1} \Pr\{P^{(-i)}_{1:m-1} \le \alpha_1, \ldots, P^{(-i)}_{r-1:m-1} \le \alpha_{r-1}, \\ &\qquad P^{(-i)}_{r:m-1} > \alpha_r, P_i > \alpha_r\} \end{aligned} \tag{3.5}$$

$$\leq \alpha \sum_{i=1}^{m} \sum_{r=1}^{m} \frac{1}{m-r+1} \Pr\{P_{1:m-1}^{(-i)} \leq \alpha_1, \ldots, P_{r-1:m-1}^{(-i)} \leq \alpha_{r-1}, P_{r:m-1}^{(-i)} > \alpha_r, P_i > \alpha_r\}$$

$$= \alpha \sum_{r=1}^{m} \Pr\{P_{1:m} \leq \alpha_1, \ldots, P_{r-1:m} \leq \alpha_{r-1}, P_{r:m} > \alpha_r\}$$

$$= \alpha[1 - \Pr\{P_{1:m} \leq \alpha_1, \ldots, P_{m:m} \leq \alpha_m\}]$$

[with $P_{m:m-1}^{(-i)} = 1$ in the second line]. This proves the theorem.

REMARK 3.1. Clearly, the step-down procedure with the constants $\alpha_i = iq/(m+\beta_m - i(1-q)), i = 1, \ldots, m$, for any $\beta_m \geq 1$, will also control the FDR under independence. For $\beta_m = 1$, these constants are the largest, as well as the power of the corresponding procedure. Therefore in our simulation study, to be discussed in the next section, we have used this sequence of constants with $\beta_m = 1$.

## 4. Performance.

We conducted a simulation study to investigate numerically the FDR control level and the power of adaptive procedures with *proven* control of the FDR. The procedures compared are the multiple-stage step-down procedure, the two-stage BH procedure, the modified Storey's procedure in [10], which we call the adaptive $\lambda$-based step-up procedure (hereafter referred to as the STS procedure), and two additional procedures that serve as benchmarks for comparing performance. These two additional procedures are the BH procedure at level $q$ and the oracle procedure that controls the FDR at the exact level $q$ but is not implementable in practice as $m_0$ is unknown.

4.1. *The configurations investigated.* The number of tests $m$ was set at $m = 64$, 512 and 4096, with the fraction of the true null hypotheses $(m_0/m)$ at 0%, 25%, 50%, 75% and 100%. The $p$-values were generated in the following way. With $Z_1, \ldots, Z_{m+1}$ distributed independently and identically as $N(0,1)$, we first let $Y_i = \sqrt{\rho} Z_{m+1} + \sqrt{1-\rho} Z_i + \mu_i$, $i = 1, \ldots, m$, and then defined the corresponding $p$-values as $p_i = 1 - \Phi(Y_i)$ , $i = 1, \ldots, m$, where $\Phi$ is the cdf of $N(0,1)$. The values of $\mu_i$ were either 0 or an equal proportion of 1, 2, 3 and 4; the values of $\rho$ were 0, 0.2, 0.5 and 0.8. The simulation results are based on 5000 replications.

4.2. *Independent test statistics.* The results of the FDR control under independence for all the procedures are given in Table 1. As expected from the theory, the FDR level of the BH procedure is exactly $qm_0/m$. The FDR level of all adaptive procedures is less conservative. In fact it is hoped that

TABLE 1
*Estimated FDR values with $\rho = 0$*

| | $m_0/m = 0.25$ | | | $m_0/m = 0.5$ | | | $m_0/m = 0.75$ | | | $m_0/m = 1$ | | |
|---|---|---|---|---|---|---|---|---|---|---|---|---|
| $m$ | 64 | 512 | 4096 | 64 | 512 | 4096 | 64 | 512 | 4096 | 64 | 512 | 4096 |
| BH | 0.013 | 0.013 | 0.013 | 0.025 | 0.025 | 0.025 | 0.036 | 0.038 | 0.037 | 0.049 | 0.048 | 0.048 |
| TS | 0.023 | 0.023 | 0.022 | 0.034 | 0.034 | 0.034 | 0.039 | 0.041 | 0.041 | 0.045 | 0.045 | 0.045 |
| MS | 0.026 | 0.026 | 0.026 | 0.036 | 0.036 | 0.036 | 0.040 | 0.043 | 0.043 | 0.048 | 0.047 | 0.047 |
| STS | 0.040 | 0.039 | 0.039 | 0.046 | 0.046 | 0.046 | 0.047 | 0.049 | 0.049 | 0.049 | 0.048 | 0.048 |
| ORC | 0.051 | 0.050 | 0.050 | 0.050 | 0.050 | 0.050 | 0.048 | 0.050 | 0.050 | 0.049 | 0.048 | 0.048 |

The standard error of estimated FDR $< 0.003$ when $m_0/m = 1$. BH is the linear step-up procedure, TS is the two-stage procedure in BKY, MS is the adaptive multiple-stage step-down procedure, STS is the modified Storey's procedure [10] and ORC is the BH procedure with the true $m_0$.

the FDR level of all adaptive procedures will be close to that of the oracle procedure and will not depend on $m_0/m$. In spite of the fact that this goal is not entirely achieved, we see that the dependence of the FDR level on $m_0/m$ is moderated and no longer linear.

Power comparisons are presented in Table 2. The power of each procedure is divided by the power of the oracle to yield relative efficiency. For independent test statistics the STS procedure is the most powerful in all simulated situations. When all hypotheses are false and $m = 4096$ the power of the multiple-stage procedure equals that of the STS procedure. In all other cases the power of the multiple-stage procedure takes the second place, sometimes very close to that of the STS procedure.

4.3. *Positively dependent test statistics.* Figure 1 presents the FDR level of the studied procedures for $\rho = 0.2, 0.5$ and 0.8. Obviously the oracle and the BH procedures do control the FDR, even though now at a lower level

TABLE 2
*Power relative to the oracle procedure for $\rho = 0$*

| | $m_0/m = 0$ | | | $m_0/m = 0.25$ | | | $m_0/m = 0.5$ | | | $m_0/m = 0.75$ | | |
|---|---|---|---|---|---|---|---|---|---|---|---|---|
| $m$ | 64 | 512 | 4096 | 64 | 512 | 4096 | 64 | 512 | 4096 | 64 | 512 | 4096 |
| BH | 0.658 | 0.656 | 0.656 | 0.776 | 0.775 | 0.776 | 0.874 | 0.873 | 0.872 | 0.941 | 0.941 | 0.941 |
| TS | 0.784 | 0.784 | 0.784 | 0.865 | 0.865 | 0.866 | 0.924 | 0.926 | 0.926 | 0.958 | 0.959 | 0.959 |
| MS | 0.871 | 0.918 | 0.953 | 0.887 | 0.890 | 0.891 | 0.933 | 0.937 | 0.938 | 0.961 | 0.968 | 0.969 |
| STS | 0.923 | 0.946 | 0.953 | 0.948 | 0.957 | 0.958 | 0.975 | 0.982 | 0.983 | 0.988 | 0.993 | 0.994 |

BH is the linear step-up procedure, TS is the two-stage procedure in BKY, MS is the adaptive multiple-stage step-down procedure, STS is the modified Storey's procedure [10] and ORC is the BH procedure with the true $m_0$.

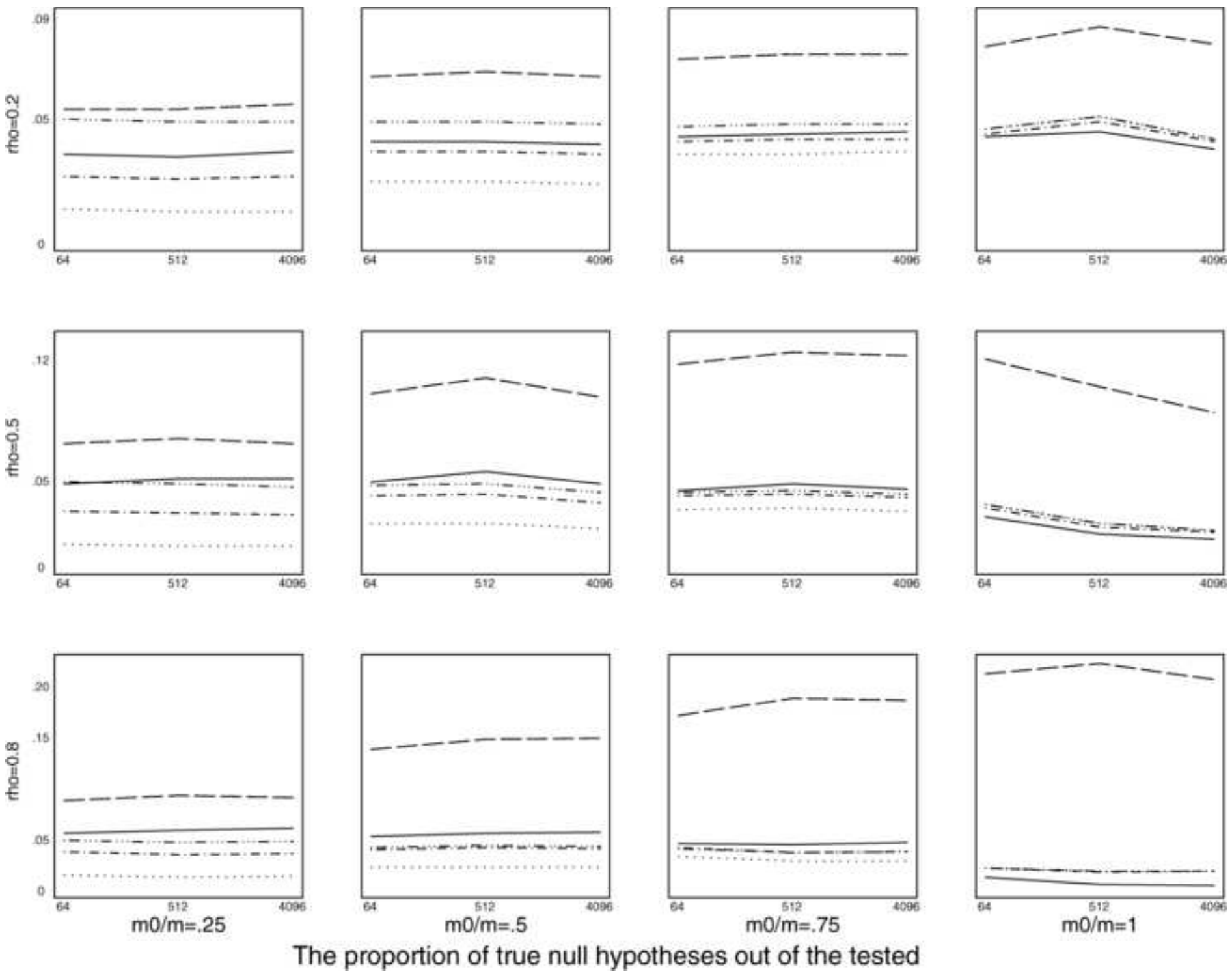

FIG. 1. *Estimated FDR values for $m = 64$, 512, 4096. Legend: MS—solid line; BH—-dotted line; TS—dotted dashed; STS—dashed; ORC—dashed triple dotted.*

as follows from theory in [5, 8]. The two-stage procedure controls the FDR below the desired level $q$ as well.

For $m_0/m = 0.25$ and 0.5 the multiple-stage procedure exceeds by little the FDR level for larger values of $\rho$. For other fractions of true null hypotheses the FDR level of multiple-stage procedure is below the desired $q$ even for higher correlations. Interestingly, for positively dependent test statistics the FDR level decreases as $m_0/m$ nears 1, in contrast to the independent case. In order to examine more closely the maximal FDR level as a function of $m_0/m$, we performed additional simulations for $m = 512$ and $\rho = 0.8$ that are presented in Figure 2. Note that the profile of FDR control as a function of $m_0/m$ of the multiple-stage step-down procedure is similar to that of the oracle procedure. It can be seen that the multiple-stage procedure is quite flat as a function of $m_0/m$ between 0.25 and 0.5 (and gets its maximal FDR level when $m_0/m \approx 0.35$). The maximal FDR level of multiple-stage procedure in our simulation is 0.061 with standard error 0.002. Therefore the FDR levels presented in Figure 1 should be close to the maximum.

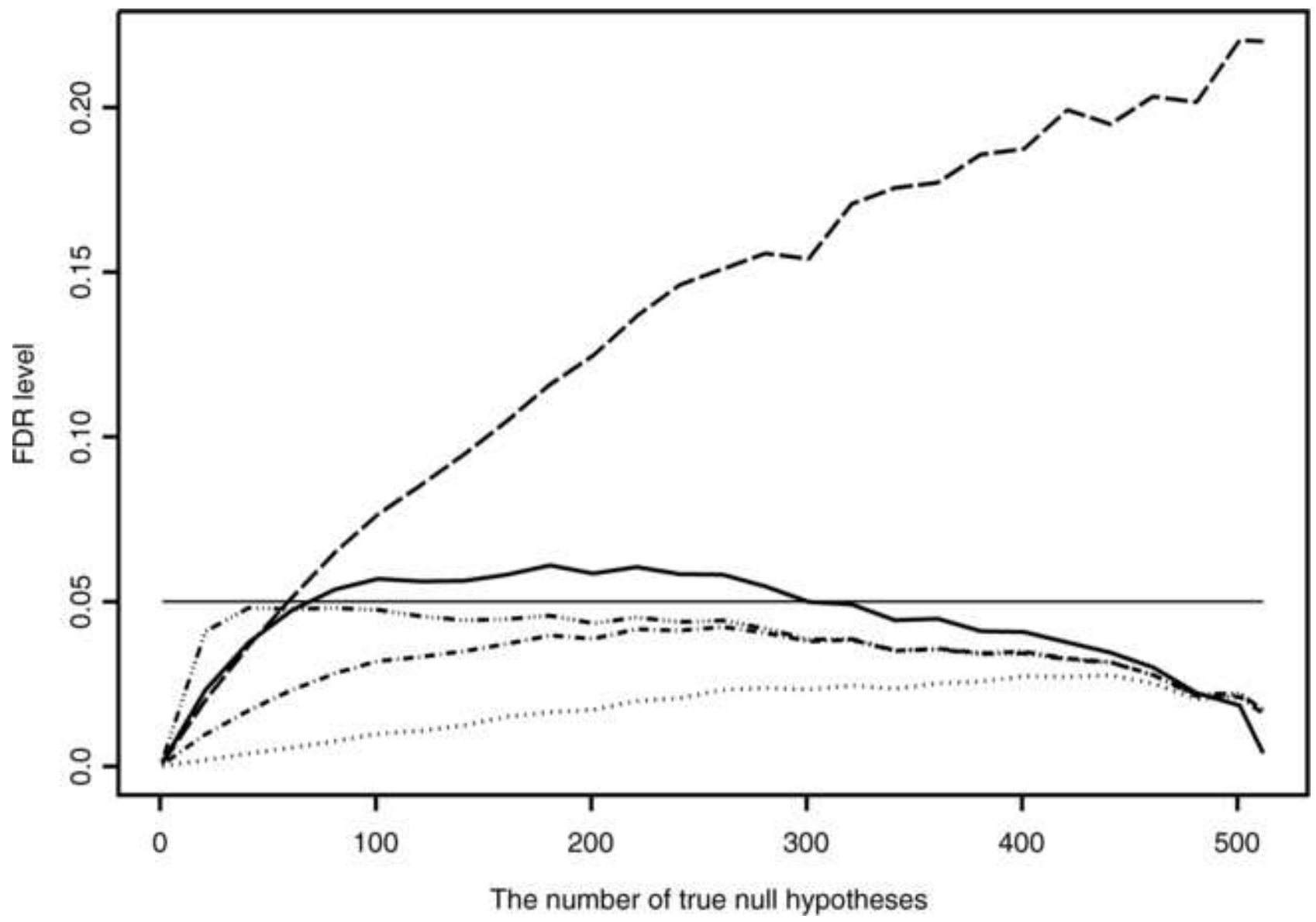

FIG. 2. *Estimated FDR values for $m = 512$, $\rho = 0.8$. Legend: MS—solid line; BH—dotted line; TS—dotted dashed; STS—dashed; ORC—dashed triple dotted.*

The FDR level of the STS procedure in extreme cases is more than four times the declared $q$. Moreover, as the results indicate, it does not seem that the FDR of the STS procedure gets closer to the nominal level as $m$ increases.

**5. Discussion.** The multiple-stage adaptive step-down procedure shows good performances among the adaptive procedures with proven FDR control so far known. It is second to the STS procedure in terms of power for independent test statistics. Moreover, its FDR control does not break down by much from the desired level for positively dependent test statistics. From our simulation study it appears that, although the FDR level of the multiple-stage procedure is a little higher than the desired level $q$ for strongly correlated test statistics, it stays below $q$ for small and moderate correlations.

The interesting question is: For positively dependent test statistics, is there some bound on FDR that can then be used, or can the constants be changed to $iq/(m + \beta_m - i(1 - q))$, as indicated in [7], with a $\beta_m > 1$ that is either a constant or a function of $m$, so that the modified procedure will control the FDR at level $q$? The answer is yes. As shown in the Appendix, if we choose $\beta_m = m(1 - q)$, the FDR is controlled for test statistics that are positively regression-dependent on subsets (PRDS) [5, 8]. Unfortunately, this procedure is much more conservative, and as $\beta_m/m$ does not converge to 0 as $m$ increases, the asymptotic optimality is lost.

An important application of the testing procedure studied here is using it as a penalized variable selection procedure with an adaptive penalty per each variable that enters the model. Unlike STS and other adaptive methods that require all $p$-values to be at hand right from the beginning, the multiple-stage adaptive procedure can be applied on the forward path of nested models, using the penalty of the form $\sigma^2 \sum_{i=1}^{k} z_{1-\alpha_i/2}^2$, and stopping at the first local minimum. In our previous simulations we found that the multiple-stage adaptive step-down procedure performs well as a variable selection procedure across a wide range of situations (even when the explanatory variables are correlated). Moreover, this procedure with the FDR level of 0.05 shows the minimax performance among the studied procedures and situations. Both empirical results quest for some theoretical justifications that are all the more plausible because of the proven results of this paper.

## APPENDIX: POSITIVE DEPENDENCE CASE

Assume that the $P_i$'s are PRDS on the subset of null $p$-values; that is,

$$E\{\phi(P_1,\dots,P_n) \mid P_i = u\} \text{ is increasing (or decreasing) in } u \tag{A.1}$$

for each $i \in I_0$ and any coordinatewise increasing (or decreasing) function $\phi$ of $(P_1,\dots,P_m)$.

THEOREM A.1. *When the $P_i$'s are PRDS, the step-down procedure with the critical values $\alpha_1 \le \cdots \le \alpha_m$ satisfying $\alpha_i/(1-\alpha_i) = iq/(m-i+\beta_m)$, $i = 1,\dots,m$, controls the FDR at $q$, for any $\beta_m \ge m(1-q)$.*

PROOF. As argued following (3.1) in proving Theorem 1.1, we first note that

$$\mathrm{FDR} \le \sum_{i\in I_0}\sum_{r=1}^{m} E\Bigg[\Pr\{P_{1:m-1}^{(-i)} \le \alpha_1,\dots,P_{r-1:m-1}^{(-i)} \le \alpha_{r-1} \mid P_i\} \times \left\{\frac{I(P_i \le \alpha_r)}{r} - \frac{I(P_i \le \alpha_{r-1})}{r-1}\right\}\Bigg].$$

Let

$$\psi(P_i) = \Pr\{P_{1:m-1}^{(-i)} \le \alpha_1,\dots,P_{r-1:m-1}^{(-i)} \le \alpha_{r-1} \mid P_i\},$$

so that, for each $i \in I_0$ and $1 \le r \le m$, the expectation under the double summation in the right-hand side of (A.2) can be expressed as follows:

$$E\left[\psi(P_i) I(P_i \le \alpha_r)\left\{\frac{1}{r} - \frac{I(P_i \le \alpha_{r-1})}{r-1}\right\}\right]. \tag{A.2}$$

Now, note that $\psi(P_i)$ is decreasing, whereas $\frac{1}{r} - \frac{I(P_i \le \alpha_{r-1})}{r-1}$ is increasing in $P_i$. Therefore, the expectation in (A.2) is less than or equal to

$$\begin{aligned}
&\frac{E\{\psi(P_i)I(P_i \le \alpha_r)\}}{\alpha_r} E\left[\left\{\frac{1}{r} - \frac{I(P_i \le \alpha_{r-1})}{r-1}\right\} I(P_i \le \alpha_r)\right] \\
&\quad = \frac{E\{\psi(P_i)I(P_i \le \alpha_r)\}}{\alpha_r}\left\{\frac{\alpha_r}{r} - \frac{\alpha_{r-1}}{r-1}\right\} \qquad \text{(A.3)} \\
&\quad \le \frac{\alpha_r}{r} - \frac{\alpha_{r-1}}{r-1}
\end{aligned}$$

as $\alpha_r/r$ is increasing in $r$. Thus, with $m_0$ true nulls, we have

$$\text{FDR} \le m_0 \sum_{r=1}^{m}\left\{\frac{\alpha_r}{r} - \frac{\alpha_{r-1}}{r-1}\right\} = \frac{m_0 \alpha_m}{m} = \frac{m_0 q}{\beta_m + mq} \le q, \qquad \text{(A.4)}$$

if $\beta_m \ge m(1-q)$. □

Y. Gavrilov
Y. Benjamini
Department of Statistics
and Operations Research
Tel Aviv University
Tel Aviv 69978
Israel
E-mail: gyulia@post.tau.ac.il
ybenja@post.tau.ac.il

S. K. Sarkar
Department of Statistics
Fox School of Business and Management
Temple University
Philadelphia, Pennsylvania 19122
USA
E-mail: sanat@temple.edu